\documentstyle{article}

\newtheorem{pro}{Proposition}
\newtheorem{theo}{Theorem}
\newtheorem{cor}{Corollary}

\def\ra{\longrightarrow}

\newcommand{\bc}[2]{
\left(
\begin{array}{c}{#1}\\{#2}\end{array}
\right)}

\title{
On the Cohomology of Theta Divisors of Hyperelliptic Jacobians}

\author{
Atsushi Nakayashiki\thanks{Faculty of Mathematics,
Kyushu University,
Ropponmatsu 4-2-1, Fukuoka 810-8560, Japan,
\quad
e-mail: atsushi@rc.kyushu-u.ac.jp}
}

\date{}
\begin{document}
\maketitle
\begin{abstract}
We prove that $k$-th singular cohomology group of the 
complement of the theta divisor in a 
hyperelliptic Jacobian is isomorphic to the $k$-th
fundamental representation of the symplectic group $Sp(2g,\bf{C})$.
This is one of the conjectures in the paper \cite{NS}.
\end{abstract}

\section{Introduction}
In our previous paper \cite{NS} 
the integrable system associated with the family of 
hyperelliptic curves is studied. 
The space of functions on the phase space modulo the action of
commuting integrals of motion is interpreted as some cohomology group.
Making extensive use of the character method several conjectures
are proposed on the structure of certain cohomology groups. 
Among them is the following conjecture on the 
dimension of the singular cohomology groups of the complement of the theta
divisor in a hyperelliptic Jacobian variety: 
$$
\dim\, H^k(J(X)-\Theta,{\bf C})=
\bc{2g}{k}-\bc{2g}{k-2},
\quad
k\leq g,
$$
where $g$ is the genus of the hyperelliptic curve $X$.
In this note we shall prove this formula.
More precisely the singular cohomology group
$H^k(J(X)-\Theta,{\bf C})$ is shown to be isomorphic to the
$k$-th fundamental representation of the symplectic group $Sp(2g,{\bf C})$
(see Corollary \ref{conj}) which is also conjectured in \cite{NS}.

In section 2 we determine the topological Euler characteristic
of the theta divisor of a hyperelliptic Jacobian.
Then we study the singular cohomology groups of the theta divisor and
the complement of the theta divisor in section 3.
We use the result of Bressler and Brylinski \cite{BB} which says 
that the singular cohomology groups and the intersection cohomology groups 
of the theta divisor of a hyperelliptic Jacobian are isomorphic.
For the intersection cohomology groups Poincar\'e duality holds.
Since lower cohomology gropus of the theta divisor are known from the 
Lefschetz theorem on hyperplane sections, other cohomology groups 
except that of the middle dimension are calculated using the 
Poincar\'e duality.
Finally the cohomology group of half degree is determined using the
Euler characteristic calculated in section 2.
The calculation of cohomology groups of the complement of the theta
divisor needs a bit more arguments.
\vskip3mm

\noindent
{\bf Acknowledgement.} I would like to thank K. Cho for stimulating
discussions and the interest on this work. A part of this work was done 
while the author stays at CRM in Universit\'e de Montreal
for concentration period on quantum algebras and integrability. 
I would like to thank
to this instiution for financial support and kind hospitality.

\section{Topological Euler Characteristic of Theta Divisor}
In this section we shall prove

\begin{theo}\label{Eulerchar} Let $X$ be a hyperelliptic curve 
of genus $g\geq 1$, 
$J(X)$ the Jacobian variety of $X$ and $\Theta$ the theta divisor. 
Then, the topological Euler characteristic of $\Theta$ is given by
\begin{eqnarray}
&&
\chi(\Theta)=(-1)^{g-1}
\Big(
\bc{2g}{g}-\bc{2g}{g-1}
\Big).
\end{eqnarray}
\end{theo}
\vskip5mm

We can assume that $X$ is determined by an equation $y^2=f(x)$ with $f(x)$
being a monic polynomial of degree $2g+1$ without multiple zeros. 
The curve $X$ has one point at infinity which we denote by $\infty$.
The hyperelliptic involution
$\sigma$ is defined by $\sigma(x,y)=(x,-y)$ 
and $\sigma(\infty)=\infty$.
Let $S^nX$ be the $n$-th symmetric products of $X$ and $\pi_n:X^n\ra S^nX$
the projection map.
For the sake of convenience we set $S^0X=\{\infty\}$ and 
$S^nX=\emptyset$ for $n<0$.
We use the bracket symbol to denote a point of $S^nX$:
$\pi_n(p_1,\cdots,p_n)=[p_1,\cdots,p_n]$.

The Jacobian variety $J(X)$ can be considered as the group 
of linear equivalence classes of divisors of degree zero.
Let $\varphi_n$ be the Abel-Jacobi map
\begin{eqnarray}
\varphi_n&:&S^nX \ra J(X)
\nonumber
\\
&&
[p_1,\cdots,p_n]\mapsto p_1+\cdots+p_n-n\infty
\nonumber
\end{eqnarray}
and $W_n$ the image of $\varphi_n$, $W_n=\varphi_n(S^nX)$ for $n\geq 1$.
We set $W_0=\varphi_0(\infty)=\{0\}$ and $W_n=\emptyset$ for $n<0$.
By Riemann's theorem $\Theta=W_{g-1}+\kappa$,
where $\kappa$ is Riemann's constant with respect to the choice of the
base point $\infty$.

For a hyperelliptic curve the structure of $W_k$ is well understood.
To see it we recall some notations on special divisors from \cite{ACGH}.
Define
$$
C^r_d=\{D\in S^dX\, \vert\, \dim H^0(X,{\cal O}(D))\geq r+1\},
\quad
W^r_d=\varphi_d(C^r_d).
$$
In particular $C^0_d=S^dX$, $W^0_d=W_d$.

A linear system ${\cal D}={\bf P}V$, $V\subset H^0(X,{\cal O}(D))$, is called
$g^r_d$ if ${\rm deg}\, D=d$ and $\dim V=r+1$. 
Here ${\bf P}V$ means the projectivization of the vector space $V$.

We refer the following fact about $g^r_d$ 
on a hyperelliptic curve (\cite{ACGH}, p13).

\vskip3mm
\begin{theo}\label{grd} Any complete $g^r_d$ on $X$ is of the form
$$
rg^1_2+p_1+\cdots+p_{d-2r},
$$
where $p_1,\cdots,p_{d-2r}\in X$ and $p_i\neq \sigma(p_j)$ 
for any $i\neq j$.
\end{theo}
\vskip3mm

The only $g^1_2$ on $X$ is the complete linear system
$$
|2\infty|=\{p+\sigma(p)\,\vert\, p\in X\}.
$$
Thus, by Theorem \ref{grd},
\begin{eqnarray}
&&
C^r_d
=
\{[p_1,\sigma(p_1),\cdots,p_r,\sigma(p_r),p_{r+1},\cdots,p_{d-r}]
\vert p_1,\cdots,p_{d-r}\in X\}.
\nonumber
\end{eqnarray}
This means, in particular, that $W^r_d=W_{d-2r}$ and 
$C^r_d=\varphi_d^{-1}(W^r_d)$.
Thus $W_{d-2}$ is exactly the singular locus of $W_d$ by
Riemann-Kempf's singularity theorem (\cite{Mum}, p56).
Notice that $C^r_d=\emptyset$ for $d-2r<0$.

We have the diagram
$$
\begin{array}{ccc}
S^kX&\stackrel{\varphi_k}{\ra}&W_k\\
\cup&\quad&\cup\\
C^1_k&\ra&W_{k-2}.
\end{array}
$$
By Theorem \ref{grd} and Riemann-Kempf's singularity theorem,
the isomorphism 
\begin{eqnarray}
&&
S^kX - C^1_k \simeq W_k - W_{k-2}
\nonumber
\end{eqnarray}
holds.
Since $S^kX$, $C^1_k$, $W_k$, $W_{k-2}$ are all 
algebraic varieties, we have ({\it cf.}\cite{F}, p95)
\begin{eqnarray}
&&
\chi(S^kX)-\chi(C^1_k)=\chi(W_k)-\chi(W_{k-2}).
\label{Fulton}
\end{eqnarray}
We set $\chi(\emptyset)=0$. Then this equation is valid for all $k\geq 0$.
We sum up (\ref{Fulton}) over $k$ and get
\begin{eqnarray}
\chi(\Theta)=\chi(W_{g-1})=
\sum_{k\geq 0}\chi(S^{g-1-2k}X)
-\sum_{k\geq 0}\chi(C^1_{g-1-2k}).
\label{kaitheta1}
\end{eqnarray}
\vskip5mm

\begin{pro}\label{c1n}
\begin{eqnarray}
&&
\chi(C^1_n)=2\chi(S^{n-2}X)-\chi(S^{n-4}X).
\label{kaitheta2}
\end{eqnarray}
\end{pro}
\vskip3mm
\noindent
{\it Proof}. Notice that 
\begin{eqnarray}
&&
C^1_2=\{[p,\sigma(p)] \vert p\in X\}
\simeq X/\{1,\sigma\}\simeq {\bf P}^1.
\nonumber
\end{eqnarray}
We identify ${\bf P}^1$ with $C^1_2$ by this isomorphism.
We define the map $\tau^{(n)}_k$ by
\begin{eqnarray}
&&
\tau^{(n)}_k: S^k{\bf P}^1\times S^{n-2k}X \ra C^k_n,
\nonumber
\\
&&
\tau^{(n)}_k
(\big[[p_1,\sigma(p_1)],\cdots,[p_k,\sigma(p_k)]\big],[p_{k+1},\cdots,p_{n-k}])
=
\nonumber
\\
&&
[p_1, \sigma(p_1), \cdots, p_k, \sigma(p_k), p_{k+1}, \cdots, p_{n-k}].
\nonumber
\end{eqnarray}
Then
\begin{eqnarray}
&&
(\tau^{(n)}_k)^{-1}(C^{k+1}_n)=
S^k{\bf P}^1\times C_{n-2k}^1,
\nonumber
\end{eqnarray}
and 
\begin{eqnarray}
&&
(S^k{\bf P}^1\times S^{n-2k}X) - (\tau^{(n)}_k)^{-1}(C^{k+1}_n)
\simeq
C^k_n - C^{k+1}_n,
\nonumber
\end{eqnarray}
by Theorem \ref{grd}.
Thus
\begin{eqnarray}
&&
\chi(S^k{\bf P}^1)\big(\chi(S^{n-2k}X)-\chi(C_{n-2k}^1)\big)=
\chi(C^k_n)-\chi(C^{k+1}_n).
\label{kaitheta3}
\end{eqnarray}
\vskip3mm

Notice that $\chi(S^k{\bf P}^1)=k+1$.
Summing up (\ref{kaitheta3}) over $k$ we have
\begin{eqnarray}
&&
\chi(C^1_n)=
\sum_{k\geq1}(k+1)
\big(
\chi(S^{n-2k}X)-\chi(C_{n-2k}^1)
\big).
\nonumber
\end{eqnarray}
The proposition can be easily proved by the induction on $n$ using
this equation.
$\Box$
\vskip3mm

It follows from (\ref{kaitheta1}), (\ref{kaitheta2}) that
\begin{eqnarray}
&&
\chi(\Theta)=\chi(S^{g-1}X)-\chi(S^{g-3}X).
\nonumber
\end{eqnarray}
The Euler characteristic of the symmetric products of curves of
$g\geq 1$ is known (\cite{Mac}):
\begin{eqnarray}
&&
\chi(S^nX)=(-1)^n\bc{2g-2}{n}.
\label{Macd}
\end{eqnarray}
The claim of the theorem follows from this.
$\Box$
\vskip2mm

\vskip3mm
\begin{cor}
Suppose that $X$ is hyperelliptic and $d\leq g-1$. Then
$$
\chi(C^r_d)=(-1)^d\Big((r+1)\bc{2g-2}{d-2r}-r\bc{2g-2}{d-2r-2}\Big).
$$
\end{cor}
\noindent
{\it Proof.} From (\ref{kaitheta3})
\begin{eqnarray}
&&
\chi(C^r_d)=\sum_{k\geq r}(k+1)
\big(
\chi(S^{d-2k}X)-\chi(C_{d-2k}^1)
\big).
\label{crdeq}
\end{eqnarray}
Substituting (\ref{kaitheta2}) to (\ref{crdeq}) we get
\begin{eqnarray}
&&
\chi(C^r_d)=(r+1)\chi(S^{d-2r}X)-r\chi(S^{d-2r-2}X).
\nonumber
\end{eqnarray}
Then the corollary follows from (\ref{Macd}).
$\Box$

\vskip3mm
\begin{cor}\label{wk}
Suppose that $X$ is hyperelliptic and $d\leq g-1$. Then
$$
\chi(W_d)=(-1)^d\Big(\bc{2g-2}{d}-\bc{2g-2}{d-2}\Big).
$$
\end{cor}

The proof of Corollary \ref{wk} is similar to that of Theorem \ref{Eulerchar}
and we leave it to the reader.

\section{Cohomologies of the theta divisor and the complement of the 
theta divisor}
In this section we assume that $(J,\Theta)$ is any principally polarized
abelian variety such that $\Theta$ is irreducible unless otherwise stated.
By the Poincar\'e-Lefschetz duality
\begin{eqnarray}
&&
H_k(J,\Theta,{\bf C})\simeq H^{2g-k}(J-\Theta,{\bf C}).
\nonumber
\end{eqnarray}
Since $J-\Theta$ is an affine algebraic variety, 
\begin{eqnarray}
&&
H^{k}(J-\Theta,{\bf C})=0, \quad k>g.
\nonumber
\end{eqnarray}
Then from the long exact sequence of homologies for the triple 
$\Theta\subset J \subset (J,\Theta)$ we have
\begin{eqnarray}
&&
H_k(\Theta,{\bf C})\simeq H_k(J,{\bf C}), \quad 0\leq k\leq g-2,
\label{equality1}
\end{eqnarray}
and the exact sequence
\begin{eqnarray}
&&
0\ra H^1(J,{\bf C}) \ra H^1(J-\Theta,{\bf C}) \ra H_{2g-2}(\Theta,{\bf C})
\ra H^2(J,{\bf C}) \ra
\label{longexact}
\\
&&
\cdots \ra H_{2g-k}(\Theta,{\bf C})\ra H^k(J,{\bf C}) 
\ra H^k(J-\Theta,{\bf C}) \ra
\cdots.
\nonumber
\end{eqnarray}

Let us study the kernel of the restriction map
\begin{eqnarray}
&&
\iota_k: H^k(J,{\bf C}) \ra H^k(J-\Theta,{\bf C}),
\nonumber
\end{eqnarray}
for $1\leq k\leq g$.
To this end we shall use the analytic description of $H^k(J,{\bf C})$.

Consider the exact sequence of sheaves:
\begin{eqnarray}
&&
0 \ra {\bf C} \ra {\cal O}_{J}(\ast \Theta) \stackrel{d}{\ra}
d\big({\cal O}_{J}(\ast \Theta)\big) \ra 0,
\label{sheafexact}
\end{eqnarray}
where ${\cal O}_{J}(\ast \Theta)$ is the sheaf of meromorphic functions
with poles only on $\Theta$ and $d$ is 
the holomorphic exterior differentiation map.
Since $H^k(J,{\cal O}_{J}(\ast \Theta))=0$ for $k>0$, we have
\begin{eqnarray}
&&
H^1(J,{\bf C})\simeq 
{
H^0(J, d\big({\cal O}_{J}(\ast \Theta)\big))
\over
dH^0(J,{\cal O}_{J}(\ast \Theta))
},
\label{secondkind}
\end{eqnarray}
by the long cohomology exact sequence of (\ref{sheafexact}).

Write $J$ as 
\begin{eqnarray}
&&
J(X)={\bf C}^g/{\bf Z}^g+\Omega{\bf Z}^g,
\nonumber
\end{eqnarray}
where $\Omega$ is a point on the Siegel upper half space, that is,
$\Omega$ is a $g\times g$ matrix with positive definite imaginary part.
We denote the coordinate system of ${\bf C}^g$ by $(z_1,\cdots,z_g)$.
Let $\theta(z)$ be Riemann's theta function
\begin{eqnarray}
&&
\theta(z)=\sum_{n\in {\bf Z}^{g}}\exp(\pi i {}^tn\Omega n+2\pi i {}^tn z).
\nonumber
\end{eqnarray}
The divisor $\Theta$ is, by definition, the zero set of $\theta(z)$.
Define $\zeta_i(z)$ by
\begin{eqnarray}
&&
\zeta_i(z)={\partial \over \partial z_i}\log\, \theta(z).
\nonumber
\end{eqnarray}

It is obvious that 
\begin{eqnarray}
&&
dz_i, \quad d\zeta_i \in
H^0(J, d\big({\cal O}_{J}(\ast \Theta)\big)).
\nonumber
\end{eqnarray}
The following proposition is easily proved by calculating
periods.
\vskip3mm

\begin{pro} The forms $dz_i$, $d\zeta_j$, $(1\leq i,j\leq g)$
form a basis of $H^1(J,{\bf C})$.
\end{pro}
\vskip3mm

This proposition is stated for hyperelliptic Jacobian in \cite{NS}.
But it is true for any principally polarized abelian variety $(J,\Theta)$
such that $\Theta$ is irreducible.

It is known that $H^k(J,{\bf C})$ is isomorphic to the k-th exterior
products of $H^1(J,{\bf C})$,
\begin{eqnarray}
&&
H^k(J,{\bf C})\simeq \wedge^k H^1(J,{\bf C}).
\nonumber
\end{eqnarray}

We set
\begin{eqnarray}
&&
\omega=\sum_{i=1}^gdz_i\wedge d\zeta_i
\in \wedge^2 H^1(J,{\bf C}).
\nonumber
\end{eqnarray}
This form is related with the fundamental class $[\Theta]$ of $\Theta$ 
in $H^2(J,{\bf C})$ (\cite{GH}, Ch.0, \S4) by
$$
[\Theta]={1\over 2\pi \sqrt{-1}}\omega.
$$
\vskip3mm

\begin{pro}\label{kernel}
$\omega\wedge H^{k-2}(J,{\bf C})\subset {\rm Ker}\, \iota_k, 
\quad 2\leq k\leq g$.
\end{pro}
\vskip3mm
\noindent
{\it Proof.} Let $\Omega^k(\ast \Theta)$ be the sheaf of meromorphic
$k$-forms on $J$ which have poles only on $\Theta$.
By the algebraic de Rham theorem ({\it cf.} \cite{GH})
\begin{eqnarray}
&&
H^k(J-\Theta,{\bf C})\simeq
{
{\rm Ker}\big(
H^0(J,\Omega^k(\ast \Theta)) \ra H^0(J,\Omega^{k+1}(\ast \Theta))
\big)
\over
dH^0(J,\Omega^{k-1}(\ast \Theta))
}.
\nonumber
\end{eqnarray}
In terms of this description of $H^k(J-\Theta,{\bf C})$ and the description
(\ref{secondkind}) of $H^1(J,{\bf C})$, the map $\iota_k$ has 
a simple meaning.

We set $dz_I=dz_{i_1}\wedge\cdots\wedge dz_{i_k}$ and $l(I)=k$ for 
$I=(i_1, \cdots, i_k)$ etc.
Then for $dz_{I_1}\wedge d\zeta_{I_2}\in H^k(J,{\bf C})$, 
which means $l(I_1)+l(I_2)=k$ in particular,
the image $\iota_k(dz_{I_1}\wedge d\zeta_{I_2})$ 
is nothing but the meromorphic differential form on 
$J$ obtained as the exterior products of the meromorphic
one forms $dz_{i_1}$, ..., $d\zeta_{j_r}$, where $I_1=(i_1,\cdots,i_k)$,
$I_2=(j_1,\cdots,j_r)$.
Since
\begin{eqnarray}
&&
\omega=\sum_{i,j=1}^g
{\partial^2 \over \partial z_i\partial z_j}
\log\theta(z) dz_i\wedge dz_j
=0
\nonumber
\end{eqnarray}
as a meromorphic differential form on $J$, the proposition is proved.
$\Box$
\vskip2mm

Since $\Theta$ is irreducible,
\begin{eqnarray}
&&
H_{2g-2}(\Theta,{\bf C})={\bf C}\mu_{\Theta},
\label{fundclass}
\end{eqnarray}
where $\mu_{\Theta}$ is the fundamental class of $\Theta$ 
(\cite{Iv}, Ch V, Theorem 3.4).
It follows from (\ref{longexact}), (\ref{fundclass}) and
Proposition \ref{kernel} that the map
$$
H_{2g-2}(\Theta,{\bf C}) \ra H^2(J,{\bf C})
$$
is injective and
\begin{eqnarray}
&&
H^1(J,{\bf C}) \simeq H^1(J-\Theta,{\bf C}).
\nonumber
\end{eqnarray}
This isomorphism is proved in \cite{NS} in another way.

\vskip3mm
\begin{cor} 
${\rm dim}\big({\rm Ker}\iota_k\big)\geq \bc{2g}{k-2}, 
\quad 2\leq k\leq g$.
\end{cor}
\vskip3mm
\noindent
{\it Proof.} Since the map 
\begin{eqnarray}
&&
\omega\wedge:\wedge^k H^1(J,{\bf C}) \ra \wedge^{k+2} H^1(J,{\bf C})
\nonumber
\end{eqnarray}
is injective for $0\leq k\leq g-1$ by the hard Lefschetz theorem, 
we have the assertion of the corollary.
$\Box$
\vskip3mm

\begin{cor} If
\begin{eqnarray}
&&
\dim\, H_{k}(\Theta)\leq \bc{2g}{k+2}
\label{inequality1}
\end{eqnarray}
for some $g\leq k\leq 2g-2$,
the map $H_{k}(\Theta,{\bf C}) \ra H^{2g-k}(J,{\bf C})$ is injective and the
equality holds in (\ref{inequality1}). 
\end{cor}
\vskip3mm

\begin{cor} If 
\begin{eqnarray}
&&
\dim\, H_{k}(\Theta)\leq \bc{2g}{k+2},
\label{inequality2}
\end{eqnarray}
for all $g\leq k\leq 2g-2$, then we have the exact sequences
\begin{eqnarray}
&&
0 \ra H_{2g-l}(\Theta,{\bf C}) \ra H^l(J,{\bf C}) \ra
H^l(J-\Theta,{\bf C}) \ra 0,\quad
l\leq g-1,
\nonumber
\\
&&
0\ra H_{g}(\Theta,{\bf C}) \ra H^g(J,{\bf C}) 
\ra H^g(J-\Theta,{\bf C}) \ra 
\nonumber
\\
&&
\ra H_{g-1}(\Theta,{\bf C})
\ra H^{g+1}(J,{\bf C}) \ra 0.
\nonumber
\end{eqnarray}
\end{cor}
\vskip5mm

\begin{theo}\label{main} Let $X$ be a hyperelliptic curve and $(J(X),\Theta)$
the Jacobian variety of $X$ and its theta divisor.
Then
$$
\dim\, H_k(\Theta,{\bf C})=
\left\{
\begin{array}{ll}
\bc{2g}{k},& \quad 0\leq k\leq g-2\\
\bc{2g}{k+2},& \quad g-1\leq k\leq 2g-2.
\end{array}\right.
$$
\end{theo}
\vskip3mm
\noindent
{\it Proof.} For $0\leq k\leq g-2$ the claim is already established
by (\ref{equality1}).
To prove the remaining part we use the following result of 
\cite{BB} (Proposition 3.2.1 (2)).
\vskip3mm

\begin{theo}(\cite{BB}) If $X$ is hyperelliptic
\begin{eqnarray}
&&
H^k(\Theta,{\bf C})\simeq IH^k(\Theta),
\nonumber
\end{eqnarray}
where $IH^k(\Theta)$ is the intersection cohomology group
of $\Theta$ with middle perversity (\cite{GM}).
\end{theo}
\vskip3mm

Since the Poincar\'e duality holds for the intersection cohomology groups 
(\cite{GM2,GM}),
so is for the singular cohomology groups of 
$\Theta$ if $X$ is hyperelliptic.
Thus
\begin{eqnarray}
&&
\dim\, H^k(\Theta,{\bf C})=
\dim\, H^{2g-2-k}(\Theta,{\bf C})=
\bc{2g}{k+2},
\quad
g\leq k\leq 2g-2.
\nonumber
\end{eqnarray}
Finally 
\begin{eqnarray}
\dim\, H^{g-1}(\Theta,{\bf C})
&=&
(-1)^{g-1}\big(
\chi(\Theta)-\sum_{k\neq g-1}(-1)^k \dim\, H^k(\Theta,{\bf C})
\big)
\nonumber
\\
&=&
\bc{2g}{g+1}.
\nonumber
\end{eqnarray}
Thus Theorem \ref{main} is proved.
$\Box$
\vskip3mm

\begin{cor}\label{conj} If $X$ is hyperelliptic,
\begin{eqnarray}
&&
H^k(J(X)-\Theta,{\bf C})\simeq 
{
\wedge^k H^1(J(X),{\bf C})
\over
\omega \wedge^{k-2} H^1(J(X),{\bf C})
},
\quad
0\leq k\leq g.
\label{fundrep}
\end{eqnarray}
In particular
\begin{eqnarray}
&&
\dim\, H^k(J(X)-\Theta,{\bf C})=
\bc{2g}{k}-\bc{2g}{k-2}.
\nonumber
\end{eqnarray}
\end{cor}
\vskip3mm

With the intersection form, $H^1(J(X),{\bf C})$ is a symplectic vector space
of dimension $2g$. It is isomorphic to the first fundamental representation
of the symplectic group $Sp(2g,{\bf C})$.
Then it is well known that the right hand side of (\ref{fundrep}) is
isomorphic to the $k$-th fundamental representation of $Sp(2g,{\bf C})$
(see \cite{FH} for example).

Corollary \ref{conj} is one of conjectures in \cite{NS}.

\end{document}